\documentclass[10pt,reqno,twoside]{article}

\usepackage{amsfonts}
\usepackage{amssymb}
\usepackage{amsbsy}
\usepackage{amsmath,cite,graphicx}
\usepackage{subcaption}
\usepackage[numbers,sort&compress]{natbib}
\usepackage[top=20mm,bottom=20mm,left=15mm,right=15mm]{geometry}
\newcommand\blfootnote[1]{%
  \begingroup
  \renewcommand\thefootnote{}\footnote{#1}%
  \addtocounter{footnote}{-1}%
  \endgroup
}

\usepackage[T1]{fontenc}
\newtheorem{theorem}{Theorem}[section]
\numberwithin{equation}{section}
\numberwithin{figure}{section} 
\newtheorem{example}[theorem]{Example}

\begin{document}

\title{A new technique for solving singular IVPs of Emden-Fowler type}

\author{Necdet B\.{I}LD\.{I}K$^{1}$\footnote{Corresponding author} \ and Sinan DEN\.{I}Z$^{1}$ }

\date{{\small$^{1}$ Department of Mathematics, Faculty of Art and Sciences, Celal Bayar University, 45040 Manisa, Turkey.\\ e-mails: n.bildik@cbu.edu.tr, sinan.deniz@cbu.edu.tr.
}}

\maketitle

\blfootnote{
\rule{0.63cm}{0cm}\emph{Mathematics Subject Classification:} 34E20,	34A12, 85A04}

\begin{abstract}
A new analytic approximate technique for addressing nonlinear problems, namely the optimal perturbation iteration method, is introduced and implemented to singular initial value Lane-Emden type problems to test the effectiveness and performance of the method. This technique provides us to adjust the convergence regions when necessary.Comparing different methods reveals that the proposed method is highly accurate and  has great potential to be a new kind of powerful analytical tool for nonlinear differential equations.
\end{abstract}

\section{Introduction}
Many problems of science and engineering lead to different types of differential equations and it is still very hard to solve them in the presence of strong nonlinearity.
In order to cope with this issue, there has been much attention devoted  to investigate better and more efficient solution techniques such as homotopy decomposition method \cite{atangana2013analytical,atangana2013solving}, auxiliary equation method \cite{pinar2015observations,pinar2015analytical},homotopy perturbation method \cite{sakar2016numerical,elbeleze2013homotopy},Taylor collocation method \cite{bildik2015comparison,deniz2014comparison,sezer1994taylor},Sumudu transform method\cite{bulut2013analytical} and Adomian decomposition method \cite{bulut2004numerical,abbasbandy2003improving,bildik2006use,bildik2006solution}.

Emden-Fowler equation is one of the most important differential equations of mathematical physics. It distinctively characterizes many scientific phenomena. The generalized Emden-Fowler equation is defined in the following form

\begin{equation}\label{Generalform}
{y}''+\alpha (x){y}'+\beta (x)\gamma (y)=0	
\end{equation}
			
\noindent where $\alpha (x),\beta (x),\gamma (y)$ are some arbitrary functions.For different $\gamma(y)$, the Eq. (\ref{Generalform}) has been subject of many studies in the literature such as the theory of stellar structure, thermionic currents and isothermal gas spheres. When  $\alpha (x)=\frac{k}{x},\beta (x)={{\beta }_{0}}{{x}^{r}},\gamma (y)={{y}^{s}}$ ($k$ and ${{\beta }_{0}}$ are constants, $s$ and $r$ are  real numbers), Eq.(\ref{Generalform}) reduces to the classic Emden-Fowler equation:

\begin{equation}\label{classicEmdenFowlerequation}
{y}''+\frac{k}{x}{y}'+{{\beta }_{0}}{{x}^{r}}{{y}^{s}}=0.	
\end{equation}
			
\noindent Furthermore, by choosing $r=0$ and $k=2$ , we get the standard Lane-Emden equation
\begin{equation}\label{standardLaneEmdenequation}
{y}''+\frac{2}{x}{y}'+{{\beta }_{0}}{{y}^{s}}=0
\end{equation}

\noindent which arises in astrophysics. Eq. (\ref{classicEmdenFowlerequation}) is also used to model the thermal behavior of a spherical cloud of gas acting under the mutual attraction of its molecules.  Many analytical techniques have been considered by various researchers to obtain the approximate solutions for these types of  equations \cite{liao2003new,he2003variational,parand2010approximation,wong1975generalized}.

Perturbation iteration method has been developed recently by Pakdemirli et al and it  has been successfully implemented to some strongly nonlinear systems\cite{aksoy2010new,aksoy2012new,csenol2013perturbation,timuccin2013new}.In this paper, new optimal perturbation iteration algorithms(OPIA) are constructed based on perturbation iteration method. Then the simplest OPIA is applied to obtain  reliable approximate solutions to the Lane-Emden equations
  \begin{equation}\label{LaneEmdenequation}
{y}''+\frac{2}{x}{y}'+\beta (x)\gamma (y)=0
\end{equation}
with different choices of  $\beta(x),\gamma(y)$. It is also proved that this method enables us to control the convergence of solution series for the given illustrations.

\section{Optimal Perturbation Iteration Algorithms}
In this section, the following formulation is given to explain the basic concept of OPIAs for second order differential equations.

\textbf{(a)} Write the governing differential equation as:
\begin{equation} \label{governingdifferential}
F({y}'',{y}',y,\varepsilon )=A(y)+g(x)=0
\end{equation}

\noindent where  $\varepsilon$  is the auxiliary perturbation parameter, $y=y(x)$ is the unknown function and $g(x)$ is the source term. Furthermore,(\ref{governingdifferential}) can be decomposed into $A(y)=L(y)+N(y)$ where $L,N$ are linear and nonlinear parts respectively. Here we have a great freedom to choose linear part $L$.

\vglue0.2cm
\textbf{(b)} Approximate solution is taken as

\begin{equation} \label{approximatesolution}
{y_{n+1}}={y_{n}}+\varepsilon{\left( {y_c} \right)}_{n}
\end{equation}

\noindent with one correction term in the perturbation expansion. Substituting the Eq. (\ref{approximatesolution}) into (\ref{governingdifferential}) and expanding Taylor series will create the algorithms. OPIAs are classified  with respect to the degrees of derivatives in the Taylor expansions($m$). Briefly, this process is represented as OPIA-$m$. OPIA-1 and OPIA-2 are constructed below.
\vglue0.2cm
\textit{\underbar{OPIA-$1$ }}\underbar{}

\vglue0.2cm
\noindent Inserting (\ref{approximatesolution}) into the nonlinear part of the (\ref{governingdifferential}) and expanding in a Taylor series  with first derivatives yields
\begin{eqnarray} \label{fOPIA1}
N({y_n}^{\prime \prime},{y_n}^{\prime},y_n,0)+N_y{(y_c)_n}\varepsilon +N_{y'}{(y_c')}_n\varepsilon+N_{y''}{(y_c'')}_n\varepsilon +N_{\varepsilon}\varepsilon=0
\end{eqnarray}
\noindent or by reorganizing
\begin{eqnarray} \label{OPIA1}
{{\left( {y_c}^{\prime \prime}\right)}_{n}}+\frac{{{N}_{{{y}'}}}}{{{N}_{{{y}''}}}}{{\left( {{y}_{c}}^{\prime } \right)}_{n}}+\frac{{{N}_{y}}}{{{N}_{{{y}''}}}}{{\left( {{y}_{c}} \right)}_{n}}=-\frac{\frac{N}{\varepsilon }+{{N}_{\varepsilon }}}{{{N}_{{{y}''}}}}.
\end{eqnarray}

\textit{\underbar{OPIA-$2$}}\underbar{}
\vglue0.2cm
\noindent Using the second order derivatives, we get
\begin{equation}
\begin{array}{l} N(y_{n} ^{''} ,y_{n} ^{'} ,y_{n} ,0)+N_{y} (y_{c} )_{n} \varepsilon+ N_{y'}(y_{c} ^{{'}})_{n} \varepsilon +N_{y''}(y_{c} ^{{'} {'} } )_{n} \varepsilon  + N_{\varepsilon } \varepsilon +\frac{1}{2} \varepsilon ^{2} N_{y''y''} (y_{c} ^{{'} {'} } )_{n}^{2}
+\frac{1}{2} \varepsilon ^{2} N_{y'y'}(y_{c} ^{{'} } )_{n}^{2} +\frac{1}{2} \varepsilon ^{2} N_{yy} (y_{c} )_{n}^{2} +  \\

\varepsilon ^{2} N_{y''y'} (y_{c} ^{{'} {'} } )_{n} (y_{c} ^{{'} } )_{n}+\varepsilon ^{2} N_{y'y} (y_{c} ^{{'} } )_{n} (y_{c} )_{n}
+\varepsilon ^{2} N_{y''y} (y_{c} ^{{'} {'} } )_{n} (y_{c} )_{n}+N_{\varepsilon y''} (y_{c} ^{{'} {'} } )_{n} \varepsilon ^{2}
+N_{\varepsilon y'} (y_{c} ^{{'} } )_{n} \varepsilon ^{2}+N_{\varepsilon y} (y_{c} )_{n} \varepsilon ^{2}
+\frac{1}{2} \varepsilon ^{2} N_{\varepsilon \varepsilon }=0
\end{array}
\end{equation}

\noindent or by rearranging we obtain

\begin{equation} \label{OPIA2}
\begin{array}{l} (y_{c} ^{''} )_{n} \left(\varepsilon N_{y''} +\varepsilon ^{2} N_{\varepsilon y''} \right)+(y_{c} ^{'} )_{n} \left(\varepsilon N_{y'} +\varepsilon ^{2} N_{\varepsilon y'} \right)+(y_{c} )_{n} \left(\varepsilon N_{y} +\varepsilon ^{2} N_{\varepsilon y} \right)+(y_{c} )_{n}^{2} \left(\frac{\varepsilon ^{2} }{2} N_{yy} \right)+
(y_{c} ^{''} )_{n}^{2} \left(\frac{\varepsilon ^{2}}{2} N_{y''y''} \right) \\

+(y_{c} ^{{'} } )_{n}^{2} \left(\frac{\varepsilon ^{2} }{2} N_{y'y'} \right)+(y_{c} )_{n} (y_{c} ^{'})_{n} \left(\varepsilon ^{2} N_{y'y} \right)+(y_{c} ^{{'} {'} } )_{n} (y_{c} ^{{'} } )_{n} \left(\varepsilon ^{2} N_{y'y''} \right)+(y_{c} ^{{'} {'} } )_{n} (y_{c} )_{n} \left(\varepsilon ^{2} N_{yy''} \right)=-N-N_{\varepsilon } \varepsilon -\frac{\varepsilon ^{2} N_{\varepsilon \varepsilon } }{2}.  \end{array}
\end{equation}

\noindent It should be emphasized that all derivatives and functions  are calculated at $\varepsilon=0$. To describe the iterative scheme, first correction terms $(y_c)_0$ can be computed  from (\ref{OPIA1}) and  (\ref{OPIA2}) by using an initial guess  $y_0$.


The Eqs. (\ref{OPIA1}) and  (\ref{OPIA2}) may seem very complicated at first, but it should not be forgotten that we use the general form of the differential equations of second order to illustrate the proposed method. Actually, most differential equations in literature contain only some of the nonlinear terms $y,y',y''$. So,the Eqs.(\ref{OPIA1}) and  (\ref{OPIA2}) reduce to some simple mathematical expressions in many cases.

\vglue0.2cm

\textbf{c)}  After obtaining $(y_c)_{0}$  for each OPIA-$m$, we use the following equation
\begin{equation} \label{formula}
{{y}_{n+1}}={{y}_{n}}+{{S}_{n}}(\varepsilon ){{\left( {{y}_{c}} \right)}_{n}}
\end{equation}

\noindent to increase the accuracy of the results and effectiveness of the method. Here $S_{n}(\varepsilon)$  is an auxiliary function which provides us to adjust the convergence. The choices of functions $S_{n}(\varepsilon)$   could be exponential, trigonometric,polynomial,etc. In this study, we select auxiliary function in the form
\begin{equation} \label{auxiliaryfunction}
{{S}_{n}}(\varepsilon )={{C}_{0}}+\varepsilon {{C}_{1}}+{{\varepsilon }^{2}}{{C}_{2}}+{{\varepsilon }^{3}}{{C}_{3}}+\cdots =\sum\limits_{i=0}^{n}{{{\varepsilon }^{i}}{{C}_{i}}}
\end{equation}

\noindent where  $C_1,C_2,\ldots $  are constants which are to be determined later. \\
Substituting $(y_c)_{0}$  into Eq.(\ref{formula}) yields:
\begin{equation} \label{EQ13}
{{y}_{1}}={{y}_{0}}+{{S}_{0}}(\varepsilon ){{\left( {{y}_{c}} \right)}_{0}}={{y}_{0}}+{{C}_{0}}{{\left( {{y}_{c}} \right)}_{0}}.
\end{equation}

\noindent By using initial and boundary conditions , the first approximate result is obtained:

\begin{equation} \label{far}
{{y}_{1}}=y(x,{{C}_{0}}).
\end{equation}

\noindent Repeating the similar steps by using the Eq.(\ref{far}) and setting $\varepsilon=1$, other approximate results are found as:

\begin{equation} \label{ars}
\begin{array}{l} {y_2(x,{{C}_{0}},{{C}_{1}})={{y}_{1}}+{{S}_{1}}(\varepsilon ){{\left( {{y}_{c}} \right)}_{0}}={{y}_{1}}+\left( {{C}_{0}}+{{C}_{1}} \right){{\left( {{y}_{c}} \right)}_{1}}} \\

 y_3(x,{{C}_{0}},{{C}_{1}},{{C}_{2}})={{y}_{2}}+\left( {{C}_{0}}+{{C}_{1}}+{{C}_{2}} \right){{\left( {{y}_{c}} \right)}_{2}} \\

 \qquad \vdots \\

 y_m(x,{{C}_{0}},\ldots ,{{C}_{m-1}})={{y}_{m-1}}+\left( {{C}_{0}}+\cdots +{{C}_{m-1}} \right){{\left( {{y}_{c}} \right)}_{m-1}}   \end{array}
\end{equation}

\textbf{d)} Substitute  $y_m$ into the Eq.(\ref{governingdifferential})  and the general problem results in the following residual:
\begin{equation} \label{residual}
R(x,{{C}_{1}},\ldots ,{{C}_{m-1}})=L\left( {{y}_{m}}(x,{{C}_{1}},\ldots ,{{C}_{m-1}}) \right)+N\left( {{y}_{m}}(x,{{C}_{1}},\ldots ,{{C}_{m-1}}) \right) +g(x).
\end{equation}
\noindent Obviously, when $R(x,{{C}_{1}},\ldots ,{{C}_{m-1}})=0$ then the approximation ${{y}_{m}}(x,{{C}_{1}},\ldots ,{{C}_{m-1}})$ will be the exact solution. Generally it doesn't occur in nonlinear equations. To obtain the optimum values of ${{C}_{1}},{{C}_{2}},\ldots $; one may use the method of least squares:

\begin{equation} \label{ls}
J({{C}_{1}},\ldots ,{{C}_{m-1}})=\int\limits_{a}^{b}{{{R}^{2}}}(x,{{C}_{1}},\ldots ,{{C}_{m-1}})dx
\end{equation}

\noindent where $a$ and $b$  are selected from the domain of the problem.
The constants ${{C}_{1}},{{C}_{2}},\ldots $  can also be defined from
\begin{equation} \label{als}
R({{x}_{1}},{{C}_{i}})=R({{x}_{2}},{{C}_{i}})=\cdots =R({{x}_{m}},{{C}_{i}})=0,i=1,2,\ldots ,m
\end{equation}

\noindent where ${{x}_{i}}\in (a,b)$. Putting these constants into the last one of the Eqs. (\ref{ars}), the approximate solution of order $m$ is determined.

\section{Applications}
In this section, we only consider OPIA-1 to solve the problems. Obtained results show that even  the simplest algorithm gives better solutions when compared with many other methods.
\begin{example}
Consider the following homogeneous Lane-Emden type equation \cite{yildirim2009solutions,yildirim2007solutions}:

\begin{equation} \label{Ex1}
{y}''+\frac{2}{x}y'-(4x^2+6)y=0,y(0)=1,y'(0)=0, 0 \leq x \leq 1.
\end{equation}

\vglue0.3cm
\noindent Rewrite the Eq. (\ref{Ex1}) as:
\begin{equation} \label{Eq1}
F(y'',{y}',y,\varepsilon )={y}''+\varepsilon\left(\frac{2}{x}y'-(4x^2+6)y\right)=L(y)+N(y)=0
\end{equation}
where $L={y}''$ and $N(y)=\varepsilon\left(\frac{2}{x}y'-(4x^2+6)y\right)$. Eq. (\ref{fOPIA1}) simplifies to:
\begin{equation} \label{Eq2}
N({{y}_{n}},0)+{N_y}({y_c})_{n}\varepsilon+{N_{y'}}{({y_{c}}')}_{n}\varepsilon+{{N}_{\varepsilon }}\varepsilon =0
\end{equation}
\noindent Using the linear terms with Eqs. (\ref{approximatesolution}), (\ref{Eq2}) and setting  $\varepsilon=1$ yields

\begin{equation} \label{firstalgorithm}
{(y_c)''}_n=-(y_n)''-\frac{2}{x}(y_n)'+(4x^2+6)y_n.
\end{equation}

\noindent One may select $y_0=1$  as a starting guess which satisfies the given initial conditions. Using the Eqs. (\ref{formula})and (\ref{firstalgorithm}),we have

\begin{equation}\label{frstapprsolution}
{{y}_{1}}=1+{C_0}\left(\frac{x^4}{3}+3 x^2\right).
\end{equation}

\noindent It should be emphasized that $y_1$  does not represent the first correction term; rather it is the approximate solution after the first iteration. By using the procedure mentioned in Section 2, one obtains the following approximate solutions:

\begin{equation}\label{EQ26}
y_2=1+{C_0}\left(\frac{x^4}{3}+3 x^2\right)+{\left(C_0+C_1\right)}\left[\frac{x^2}{630}  \left(15 C_0 x^6+294 C_0 x^4+595 C_0 x^2-5670 C_0+210 x^2+1890\right)\right]
\end{equation}

\begin{equation}\label{thrdappr}
\begin{array}{l}
y_3=1+{C_0}\left(\frac{x^4}{3}+3 x^2\right)+{\left(C_0+C_1\right)}\left[\frac{x^2}{630}  \left(15 C_0 x^6+294 C_0 x^4+595 C_0 x^2-5670 C_0+210 x^2+1890\right)\right]+\\
\begin{array}{l}
 \frac{{\left(C_0+C_1+C_2\right)}x^2}{727650} \times  \left[\begin{array}{l}

525 C_0^2 x^{10}+525 C_0 C_1 x^{10}+16247 C_0^2 x^8+16247 C_0 C_1 x^8+63195 C_0^2 x^6+34650 C_0 x^6+\\
63195 C_0 C_1 x^6+17325 C_1 x^6-1211133 C_0^2 x^4+679140 C_0 x^4-1211133 C_0 C_1 x^4\\
+339570 C_1 x^4-4419800 C_0^2 x^2+1374450 C_0 x^2-4419800 C_0 C_1 x^2+687225 C_1 x^2\\
+19646550 C_0^2-13097700 C_0+19646550 C_0 C_1-6548850 C_1+242550 x^2+2182950
\end{array}
\right]
\end{array}
\end{array}
\end{equation}

\noindent  Unknown constants can be obtained from the residual
\begin{equation}\label{residual1}
   R(x,{{C}_{0}},{{C}_{1}},{{C}_{2}})=L\left( {{y}_{3}}(x,{{C}_{0}},{{C}_{1}},{{C}_{2}}) \right)+N\left( {{y}_{3}}(x,{{C}_{0}},{{C}_{1}},{{C}_{2}}) \right)
\end{equation}

\noindent for the third order approximation. Using the Eq.  (\ref{als})  with $x=0.3,0.6,0.9$ yields
\begin{equation}\label{constants1}
C_0=0.3342343984217452,C_1=0.31859877627965916,C_2=0.20764389922289617
\end{equation}

\noindent Substituting these constants into the Eq.(\ref{thrdappr}), the approximate solution of the third order is obtained as:
\begin{equation}\label{expsolu}
\begin{array}{l}
{y_3}(x)=1.000000000101 + 1.0040975062848503x^2 + 0.48364749266573837x^4 +\\
 0.18568062044855377x^6 + 0.04172402028090389x^8 + 0.00419221270969861x^{10}\\ + 0.00013546572737070046x^{12}.
\end{array}
\end{equation}

\noindent Repeating the same steps one can get the following approximations:
\begin{equation}\label{iter4}
\begin{array}{l}
{y_4}(x)=1.000000000002862+0.9999999996360626x^2+0.5000000125635662x^4 \\
+0.1666664827590984x^6+0.0416680355705590x^8 + 0.0083276885417629x^{10}\\
+0.00140218809407533x^{12}+0.0001808480143782x^{14}+0.0000365711746672x^{16}.
\end{array}
\end{equation}

\begin{equation}\label{iter5}
\begin{array}{l}
{y_5}(x)=1.000000000000001+0.9999999999998731x^2+0.5000000000060835x^4\\
 +0.166666666533135x^6+0.04166666825594757x^8+0.00833332219438660x^{10}\\
 +0.001388937204147727x^{12} + 0.00019828100274015x^{14}+0.000025026069221009x^{16}\\
+2.5264083250622193\times 10^{-6}x^{18}+4.007880863972304\times 10^{-7}x^{20}.
\end{array}
\end{equation}
\end{example}

\noindent Exact solution of the problem is $y=e^{x^2}$ .This problem has been investigated by Ozis et al using variational iteration method(VIM) and homotopy perturbation method (HPM) \cite{yildirim2009solutions,yildirim2007solutions}.  Figures 3.1, 3.2 and Table 1 give important information on the convergence of the absolute errors of OPIAs and  other approximate solutions in literature. It is clear that the results obtained by OPIM are more accurate than those in \cite{yildirim2009solutions,yildirim2007solutions} .

\begin{table}\caption{Comparison of absolute errors for Example 1 at different orders of approximations.}
\centering
\begin{tabular}{|p{0.15in}|p{1in}|p{1in}|p{1in}|p{1in}|} \hline

\textit{x} & \multicolumn{2}{|p{2in}|}{ Errors for OPIA-1} & \multicolumn{2}{|p{2in}|}{Errors for VIM and HPM} \\ \hline

& $\left|y_{exact}-y_{4} \right|$   & $\left|y_{exact}-y_{5} \right|$ & $\left|y_{exact}-y_{4} \right|$ & $\left|y_{exact}-y_{5} \right|$  \\ \hline

0.1 & 3.08426$\times$10${}^{-13}$   & 2.22045$\times$10${}^{-16}$  &1.11022$\times$10${}^{-15}$     & 1.00128$\times$10${}^{-16}$ \\
0.2 & 3.85914$\times$10${}^{-13}$   & 2.22045$\times$10${}^{-16}$  &5.72165$\times$10${}^{-12}$     & 3.26406$\times$10${}^{-14}$   \\
0.3 & 5.62883$\times$10${}^{-13}$   & 1.00025$\times$10${}^{-17}$  &7.47710$\times$10${}^{-10}$     & 9.59788$\times$10${}^{-12}$   \\
0.4 & 9.64347$\times$10${}^{-13}$   & 2.44249$\times$10${}^{-15}$  &2.38451$\times$10${}^{-8}$      & 5.43454$\times$10${}^{-10}$   \\
0.5 & 1.95532$\times$10${}^{-12}$   & 1.50998$\times$10${}^{-14}$  &3.51584$\times$10${}^{-7}$      & 1.24994$\times$10${}^{-8}$   \\
0.6 & 4.77649$\times$10${}^{-12}$   & 1.01037$\times$10${}^{-13}$  &3.18608$\times$10${}^{-6}$      &1.62772$\times$10${}^{-7}$   \\
0.7 & 5.45375$\times$10${}^{-11}$   & 5.02709$\times$10${}^{-13}$  & 0.0000206568                   &1.43282$\times$10${}^{-6}$   \\
0.8 & 2.78031$\times$10${}^{-11}$   & 2.05902$\times$10${}^{-12}$  & 0.000104921                    &9.47740$\times$10${}^{-6}$   \\
0.9 & 3.08426$\times$10${}^{-10}$   & 6.38223$\times$10${}^{-12}$  & 0.000442699                    & 0.000050436  \\
  1 & 2.10201$\times$10${}^{-9}$    & 2.90235$\times$10${}^{-12}$  & 0.00161516                     & 0.000226273 \\ \hline
\end{tabular}
\end{table}

\begin{figure}
\begin{subfigure}{.5\textwidth}
  \centering
  \includegraphics[width=.9\linewidth]{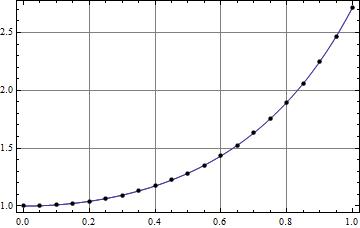}
  \caption{OPIA-1 solution of second order.($\bullet$)}
\end{subfigure}%
\begin{subfigure}{.5\textwidth}
  \centering
  \includegraphics[width=.9\linewidth]{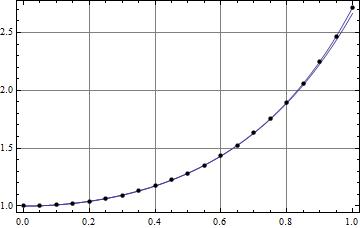}
  \caption{VIM solution of second order.($\bullet$)}
\end{subfigure}
\caption{Comparison with the exact solution(--) for
Example 1.}
\end{figure}

\begin{figure}
\begin{subfigure}{.5\textwidth}
  \centering
  \includegraphics[width=.9\linewidth]{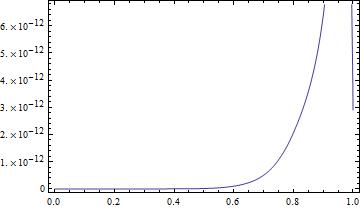}
  \caption{OPIA-1 solution of fifth order.}
\end{subfigure}%
\begin{subfigure}{.5\textwidth}
  \centering
  \includegraphics[width=.9\linewidth]{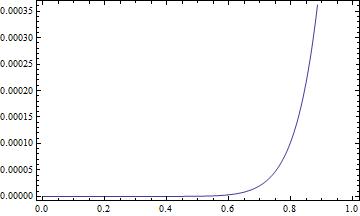}
  \caption{VIM solution of fifth order.}
\end{subfigure}
\caption{Absolute errors for OPIA-1  and the variational iteration method for Example 1.}
\end{figure}

\begin{example}
Consider the Lane-Emden type equation  \cite{khan2012solving,parand2010approximation,dehghan2008approximate}:

\begin{equation} \label{Ex2}
{y}''+\frac{2}{x}y'+e^y=0,y(0)=y'(0)=0
\end{equation}
which represents the isothermal gas spheres equation in the case that the temperature stays constant.
\vglue0.3cm
\noindent Reconsider the Eq. (\ref{Ex2}) as:
\begin{equation} \label{REx2}
F(y'',y',y,\varepsilon )={y}''+(\frac{2\varepsilon}{x}y'+e^{\varepsilon y})=L(y)+N(y)=0
\end{equation}
\noindent where $L={y}''$ and $N(y)=(\frac{2\varepsilon}{x}y'+e^{\varepsilon y})$. By using the Eqs. (\ref{approximatesolution}), (\ref{Eq2}),(\ref{REx2}) and setting  $\varepsilon=1$, we have

\begin{equation} \label{Ex2algorithm}
{(y_c)''}_n=-(y_n)''-\frac{2}{x}(y_n)'-y_n-1.
\end{equation}

\noindent By choosing $y_0=0$  as a starting guess, using the Eqs. (\ref{formula}), (\ref{Ex2algorithm}) with initial conditions, we get

\begin{equation}\label{Ex2as1}
y_1=\frac{-C_0}{2} (x^2)
\end{equation}

\begin{equation}\label{Ex2as2}
y_2=\frac{-C_0}{2} (x^2)+\frac{1}{24} \left(C_0+C_1\right) x^2 \left(C_0 x^2+36 C_0-12\right)
\end{equation}

\begin{equation}\label{Ex2as3}
\begin{array}{l}
y_3=\displaystyle\frac{-C_0}{2} (x^2)+\frac{1}{24} \left(C_0+C_1\right) x^2 \left(C_0 x^2+36 C_0-12\right)+\\
\displaystyle\frac{C_0+C_1+C_2}{720} \left[\begin{array}{l}-C_0^2 x^6-C_0 C_1 x^6-140 C_0^2 x^4+60 C_0 x^4-140 C_0 C_1 x^4+30 C_1 x^4\\
-3240 C_0^2 x^2+2160 C_0 x^2-3240 C_0 C_1 x^2+1080 C_1 x^2-360 x^2\end{array} \right]
\end{array}
\end{equation}

\begin{equation}\label{Ex2as4}
\begin{array}{l}
y_4=\displaystyle\frac{-C_0}{2} (x^2)+\frac{1}{24} \left(C_0+C_1\right) x^2 \left(C_0 x^2+36 C_0-12\right)+\\
\displaystyle\frac{C_0+C_1+C_2}{720} \left[\begin{array}{l}-C_0^2 x^6-C_0 C_1 x^6-140 C_0^2 x^4+60 C_0 x^4-140 C_0 C_1 x^4+30 C_1 x^4\\
-3240 C_0^2 x^2+2160 C_0 x^2-3240 C_0 C_1 x^2+1080 C_1 x^2-360 x^2\end{array} \right]+\\[0.5cm]

\displaystyle\frac{C_0+C_1+C_2+C_3}{604800} \left[\begin{array}{l}
15 C_0^3 x^8+15 C_0 C_1^2 x^8+30 C_0^2 C_1 x^8+15 C_0^2 C_2 x^8+15 C_0 C_1 C_2 x^8+5096 C_0^3 x^6-2520 C_0^2 x^6+\\
5096 C_0 C_1^2 x^6-840 C_1^2 x^6+10192 C_0^2 C_1 x^6-3360 C_0 C_1 x^6+5096 C_0^2 C_2 x^6-1680 C_0 C_2 x^6\\
+5096 C_0 C_1 C_2 x^6-840 C_1 C_2 x^6+422800 C_0^3 x^4-352800 C_0^2 x^4+422800 C_0 C_1^2 x^4-117600 C_1^2 x^4\\
+75600 C_0 x^4+845600 C_0^2 C_1 x^4-470400 C_0 C_1 x^4+50400 C_1 x^4+422800 C_0^2 C_2 x^4-2352C_2 x^4+\\
422800 C_0 C_1 C_2 x^4-117600 C_1 C_2 x^4+25200 C_2 x^4+8164800 C_0^3 x^2-8164800 C_0^2 x^2+816C_0 C_1^2 x^2\\
-2721600 C_1^2 x^2+2721600 C_0 x^2+16329600 C_0^2 C_1 x^2-10886400 C_0 C_1 x^2+1814400 C_1 x^2\\
-5443200 C_0 C_2 x^2+8164800 C_0 C_1 C_2 x^2-2721600 C_1 C_2 x^2+907200 C_2 x^2-302400 x^2\end{array} \right]\\
\vdots
\end{array}
\end{equation}
\end{example}

\noindent Using the Eq.  (\ref{als})  with the residual
\begin{equation}\label{residual2}
   R(x,C_0,C_1,C_2,C_3)=L\left( {{y}_{4}}(x,{{C}_{0}},{{C}_{1}},{{C}_{2}},C_3) \right)+N\left( {{y}_{3}}(x,{{C}_{0}},{{C}_{1}},{{C}_{2}},C_3) \right)
\end{equation}
\noindent the unknown constants are obtained as
\begin{equation}\label{Ex2constants1}
C_0=2.0203622551,C_1=-1.0201147822,C_2=-0.9963202221,C_3=0.020789994
\end{equation}
\noindent for the fourth order approximation.
\noindent Substituting  the Eq.(\ref{Ex2constants1}) into (\ref{Ex2as4}) yields:
\begin{equation}\label{Ex2expsolu}
\begin{array}{l}
{y_4}(x)=-0.15989962328x^2 + 0.007920058473x^4 + -0.005608897215631x^6 + 0.000039055217456x^8.
\end{array}
\end{equation}

\noindent Figure 3.3 is sketched for comparison of the higher order approximate solutions of OPIM and analytical results given by other methods in literature \cite{khan2012solving,wazwaz2011comparison}.It is also clear from the Figure 3.3 that OPIA-1 solutions are valid in larger region.

\begin{figure}
  \centering
  \includegraphics[width=.8\linewidth]{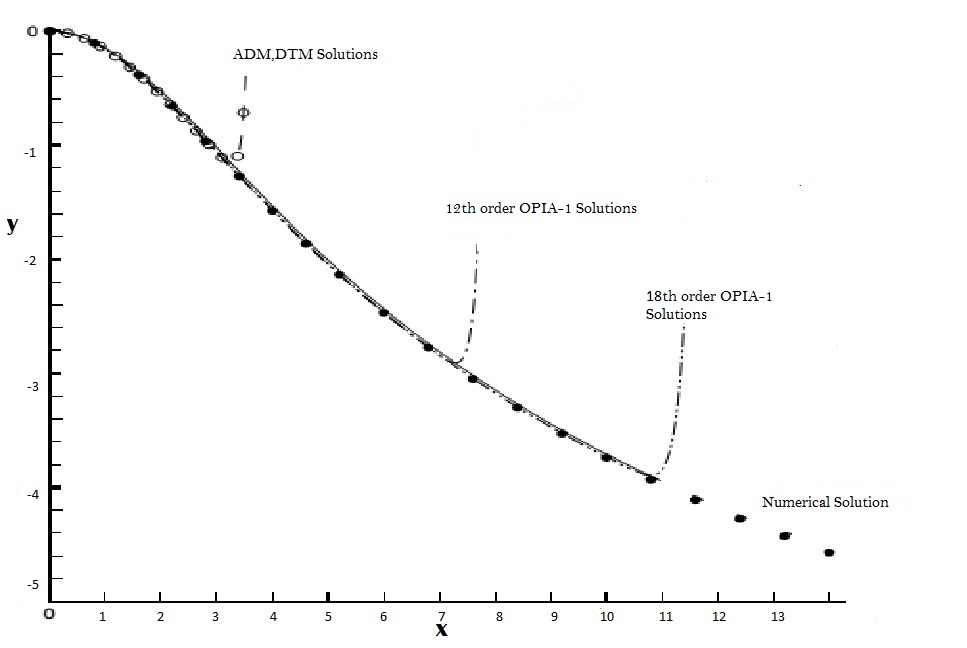}
\caption{Comparison between the results obtained by OPIM,ADM,DTM  and the numerical results for
Example 2.}
\end{figure}

\section{Conclusion}
In this paper, OPIM is applied for the first time to investigate the new approximate solutions for Lane-Emden type differential equations. This new technique provides us to optimally control the convergence of solution series. Also,it gives a very good approximation even in a few terms to these kinds of nonlinear equations. The results obtained in this paper proves that the OPIM is a very effective technique for differential equations. It is worth mentioning that, a symbolic program is necessary for successive calculations after a few iterations.Mathematica 9 has been used to overcome the complicated calculations for this present research.
\vglue1cm


\bibliographystyle{plain}
\bibliography{Bibfile}{}

\end{document}